\documentclass[12pt]{amsart}
\usepackage{amsfonts, amssymb, latexsym}
\usepackage{amssymb,amsfonts,amsmath,latexsym,cite,color, psfrag,graphicx,ifpdf,ytableau,pgfkeys,pgfopts,xcolor}
\usepackage{tikz,url}

\newtheorem*{claim*}{Claim}

\newtheorem{Main Conjecture}[theorem]{Main Conjecture}

\theoremstyle{remark}

\theoremstyle{plain}

\newcommand{\excise}[1]{}%{$\star$\textsc{#1}$\star$}

  \newcommand{\cellsize}{14}
\newlength{\cellsz} 
\newsavebox{\cell}
\sbox{\cell}{\begin{picture}(\cellsize,\cellsize)
\put(0,0){\line(1,0){\cellsize}}
\put(0,0){\line(0,1){\cellsize}}
\put(\cellsize,0){\line(0,1){\cellsize}}
\put(0,\cellsize){\line(1,0){\cellsize}}
\end{picture}}
\newcommand\cellify[1]{\def\thearg{#1}\def\nothing{}%
\ifx\thearg\nothing
\vrule width0pt height\cellsz depth0pt\else
\hbox to 0pt{\usebox{\cell} \hss}\fi%
\vbox to \cellsz{
\vss
\hbox to \cellsz{\hss$#1$\hss}
\vss}}
\newcommand\tableau[1]{\vtop{\let\\\cr
\baselineskip -16000pt \lineskiplimit 16000pt \lineskip 0pt
\ialign{&\cellify{##}\cr#1\crcr}}}
  
  \definecolor{dark-gray}{gray}{0.35}
\newcommand\yd{\tikz\draw[black,fill=yellow] (0,0) circle (.8ex);}
\newcommand\rd{\tikz\draw[red,fill=red] (0,0) circle (.8ex);}
\newcommand\dd{\tikz\draw[black,fill=dark-gray] (0,0) circle (.8ex);}
\newcommand\whd{\tikz\draw[black, fill=white] (0,0) circle (.8ex);}
\begin{document}
\pagestyle{plain}

\mbox{}
\title{An estimation method for game complexity}
\author{Alexander Yong}
\address{Dept. of Mathematics, University of Illinois at Urbana-Champaign, Urbana IL, 61801}
\email{ayong@illinois.edu}

\author{David Yong}
\address{Next Generation School, Champaign IL, 61821}
\email{david.yong124@gmail.com}

% Abstract
% We looked at a method for estimating the complexity measure of game tree size (the number of legal games). It seems effective 
% for a number of children's games such as Tic-Tac-Toe, Connect Four and Othello.  

\date{January 29, 2019}

\maketitle

We looked at a method for estimating the complexity measure of \emph{game tree size}. It seems effective for a number of 
children's games such as \emph{Tic-Tac-Toe}, \emph{Connect Four}, and \emph{Othello}.  

G.~H.~Hardy \cite[pg.~17]{hardy} estimated the game tree size (number of legal games) of Chess to be $10^{10^{50}}$ or ``in any case a second order exponential'', but gave no reasoning. Claude Shannon wrote a seminal paper \cite{Shannon} on computer Chess. Based on master games collected by the psychologist Adriaan de Groot, he estimated 
$\approx 10^3$ options per (white, black) move pair, and that an average game is $\approx 80$ \emph{plies} (half-moves). Thus, he surmised that the game tree size (and game tree complexity)  is \emph{at least} $\approx 10^{120}$. This \emph{Shannon number} is often compared to the number of atoms in the observable universe 
$\approx 10^{80}$.

We take another avenue. Given a game ${\mathcal G}$ of two players $P_1$ and $P_2$, generate a game
$g$ by uniformly at random selecting each move from legal possibilities. If $c_j=c_j(g)$ is the number of options at ply $j$, and $g$
has $N$ plies, let  $X(g)=\prod_{j=1}^N c_j$.  
Independently repeat this process, producing random games
$g_1,g_2,\ldots,g_n$. The proposed estimate for the game tree size is
\begin{equation}
\label{eqn:thetrick}
{\sf gts}({\mathcal G})\approx \frac{1}{n}\sum_{i=1}^n X(g_i).
\end{equation}
As justified below, (\ref{eqn:thetrick}) is an equality in the $n\to\infty$ limit. Our thesis is that it gives 
 \emph{fairly precise} estimates for many games of pure strategy. This estimation method is straightforward to implement, parallelizable, and space efficient. It requires no sophistication, such as use of databases (e.g., of human play, nor of endgame positions). 

To illustrate, if ${\mathcal G}=$Tic-Tac-Toe, ${\sf gts}({\mathcal G})=255168$ (\url{www.se16.info/hgb/tictactoe.htm} is our
resource). One 
random game we generated is
\[\begin{array}{c|c|c}
  \ & {\rm X} & \ \\      \hline
  \ & \ & \ \\      \hline
  \ & \ & \
\end{array} \mapsto
\begin{array}{c|c|c}
  \ & {\rm X} & \ \\      \hline
  {\rm O} & \ & \ \\      \hline
  \ & \ & \
\end{array}
\mapsto
\begin{array}{c|c|c}
  {\rm X} & {\rm X} & \ \\      \hline
  {\rm O} & \ & \ \\      \hline
  \ & \ & \
\end{array}
\mapsto
\begin{array}{c|c|c}
  {\rm X} & {\rm X} & \ \\      \hline
  {\rm O} & \ & {\rm O} \\      \hline
  \ & \ & \
\end{array}\]
\[\mapsto\begin{array}{c|c|c}
  {\rm  X} & {\rm X} & \ \\      \hline
  {\rm O} & \ & {\rm O} \\      \hline
  {\rm X} & \ & \
\end{array}
\mapsto
\begin{array}{c|c|c}
  {\rm X} & {\rm X} & \ \\      \hline
  {\rm O} & {\rm X} & {\rm O} \\      \hline
  {\rm X} & {\rm O} & \
\end{array}
\mapsto
\begin{array}{c|c|c}
  {\rm X} & {\rm X} & {\rm O} \\      \hline
  {\rm O} & {\rm X} & {\rm O} \\      \hline
  {\rm X} & {\rm O} & \
\end{array}
\mapsto
\begin{array}{c|c|c}
  {\rm X} & {\rm X} & {\rm O} \\      \hline
  {\rm O} & {\rm X} & {\rm O} \\      \hline
  {\rm X} & {\rm O} & {\rm X}
\end{array}
\]
Here, $N=9$ and $(c_1,c_2,c_3,c_4,c_5,c_6,c_7,c_8,c_9)=(9,8,7,6,5,4,3,2,1)$ and $X=9!=362880$. Estimating,
using $n=2000$ and repeating for a total of $10$ trials, gives
\[255051,260562, 252352, 256586, 250916, 256457, 257380, 251800, 257448, 248988.\]
We ``guess'' that
${\sf gts}(\text{Tic-Tac-Toe}) \approx 2.55(\pm 0.04) \times 10^5$. The ``$(\pm 0.04)$'' refers to the usual \emph{standard error of the mean}. It is comforting that this agrees with the known value.

The technique is an instance of \emph{sequential importance sampling}. This Monte Carlo method can be traced back to
 Herman Kahn and Theodore Harris' article \cite{Kahn.Harris}, who credit John von Neumann.  Let $S$ be a finite set. Assign $s\in S$ probability $p_s\in (0,1]$. Define a random variable $X:S\to {\mathbb R}$
by $X(s)=1/p_s$. From the definition of expectation,
\[
{\mathbb E}[X]=\sum_{s\in S} p_s X(s)=\sum_{s\in S} p_s\times (1/p_s)=\sum_{s\in S} 1=\#S.
\]
By the law of large numbers, $\frac{1}{n}\sum_{i=1}^n X(s_i)\to \#S$.
The application to
enumeration of $S$ was popularized by the article \cite{Knuth} of Donald Knuth. He used it 
 to estimate the number of self-avoiding walks in a grid.  It has been applied, e.g.,
by Lars Rasmussen \cite{Rasmussen} to estimate permanents, and by Joseph Blitzstein and Persi Diaconis \cite{blitz} 
 to estimate
the number of graphs with a given degree sequence.  In our use, assign $p_g=1/\prod_{j=1}^N c_j$ to each leaf node 
of depth $N$ in the game tree  of ${\mathcal G}$. Now, 
while we found no literature on the efficacy of (\ref{eqn:thetrick}) 
exactly for game tree size,  there is plenty close in point, and it may have been observed without flourish.\footnote{See \cite{otherKnuth}. Also, a different
Monte Carlo technique was applied in \cite{Allis} to estimate the number of legal \emph{positions} in games. There one uses the idea of enumerating a superset
$U$ of $S$ and sampling from $U$ \emph{uniformly} at random to estimate the probability a point is in $S$. In contrast, the probability distribution we use is far from uniform. Finally, {\sf AlphaGo} and {\sf AlphaZero} use Monte Carlo methods in Go and Chess move selection \cite{Alpha}. Indeed what is described partly forms the rudiments of an {\sf AI}: simulating many games for each choice of move and pick the one that produces
the highest estimated winning percentage.}  At any rate, we hope this letter contributes some useful experience.

It's simple to modify the above to estimate all sorts of other statistics such as ${\sf agl}({\mathcal G})$, the average game length of
${\mathcal G}$. If we define 
$Y(g)=N\cdot \prod_{j=1}^N c_j$ then 
${\sf agl}({\mathcal G})={\mathbb E}[Y]/{\mathbb E}[X]$.  
Thus, we can Monte Carlo estimate ${\sf agl}({\mathcal G})$. One knows ${\sf agl}(\text{Tic-Tac-Toe})=(5\times 1440+ 6\times 5328
+7\times 47952+8\times 72576+9\times 127872)/255168 \approx 8.25$. This is accurately estimated using $2000$ trials. In addition, we can
also estimate the percentage of wins (by either player) and draws. For draws, use the random variable 
$Z(g)=\delta_g\times
\prod_{j=1}^N c_j$ where $\delta_g=1$ if $g$ is draw, and $\delta_g=0$ otherwise. Then estimate ${\mathbb E}[Z]/{\mathbb E}[X]$. For Tic-Tac-Toe,
the draw (that is, ``cat's game'') rate is $46080/255168\approx 18.1\%$. Our simulations agree closely.\footnote{Code available at \url{https://github.com/ICLUE/Gametreesize}}

Consider ${\mathcal G}=$ \emph{Connect 4}. Though commercialized by Milton Bradley (now Hasbro) in 1974, it has a longer history.
Among its alternate names is \emph{Captain's Mistress}, stemming from folklore that the game absorbed Captain James Cook (1728--1779) during his 
historic travels. The game is played on a vertical board with seven columns of height six. $P_1$ uses $\yd$ while $P_2$ uses $\rd$.  
$P_1$ moves first and chooses a column to drop their first disk into. The players alternate. At each ply, any non-full column may be
chosen. The game terminates when there are four consecutive disks of the same color in a row, column or diagonal. 

We encode an entire game with a \emph{tableau} by recording the ply at which a disk was placed. A randomly sampled game $g$ is below; it has
$X(g)=5.59\times 10^{17}$:
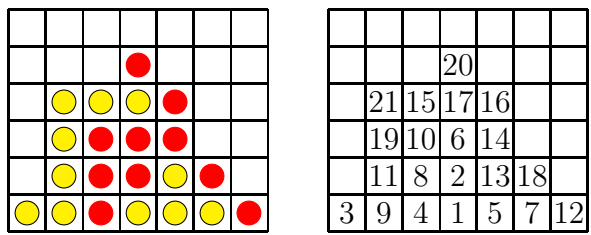
\[\tableau{ \ & \ & \ & \ & \ & \ & \ \\
 \ & \ & \ & \rd & \ & \ & \ \\
 \ & \yd & \yd & \yd & \rd & \ & \ & \\
 \ & \yd & \rd & \rd & \rd & \ & \ \\
 \ & \yd & \rd & \rd & \yd & \rd & \ \\
  \yd & \yd & \rd & \yd & \yd & \yd & \rd } \ \ \ 
\tableau{ \ & \ & \ & \ & \ & \ & \ \\
 \ & \ & \ & 20 & \ & \ & \ \\
 \ & 21 & 15 & 17 & 16 & \ & \ \\
 \ & 19 & 10 & 6 & 14 & \ & \ \\
 \ & 11 & 8 & 2 & 13 & 18 & \ \\
  3 & 9 & 4 & 1 & 5 & 7 & 12 }
\]
Thus, $\yd$ and $\rd$ correspond to odd and even labels, respectively. Since each column is increasing from bottom to top, every game of 
$N$ plies can be viewed as a distribution of $1,2,\ldots,N$
into $7$ distinguishable rooms, where each room can have at most $6$ occupants. For a fixed choice of occupancy $(o_1,\ldots,o_7)$, the number of such arrangements
is the multinomial coefficient ${N\choose o_1 \ o_2 \ \ldots \ o_7}$. Thus if $T_N$ is the total of such arrangements, then rephrasing
in terms of exponential generating series,
\[T_N = \text{coefficient of $x^N$ in \ } {N!}\left(1+\frac{x}{1!}+\frac{x^2}{2!}+\frac{x^3}{3!}+\frac{x^4}{4!}+\frac{x^5}{5!}+\frac{x^6}{6!}\right)^7.\]
Thus, 
$\#{\sf gts}({\mathcal G})\leq \sum_{N=7}^{42} T_N=40645234186579304685384521259174\approx 4.06\times 10^{31}$, 
as may be determined
quickly using a computer algebra system.

Shannon's number is an estimated lower bound for Chess' \emph{game tree complexity} ({\sf gtc}). This is
the number of leaves of the smallest full width (all nodes of each depth) decision tree determining the game-theoretic value of the initial position.
Similarly, in \cite[Section~6.3.2]{Allis}, the average game length of Connect Four \emph{in practice} is estimated to be $36$ ply with an average of $4$ legal moves/ply, whence ${\sf gtc}(\text{Connect 4})\approx  4^{36}\approx 4.72 \times 10^{21}$. 
We applied (\ref{eqn:thetrick}) with $12$ trials of the method using $n=10^8$. Based on this, 
the game tree \emph{size} appears not \emph{so} far from the upper bound:
\[{{\sf gts}(\text{Connect 4}) \approx 8.34(\pm 0.05)\times 10^{28}.}\]
%8.396066895967572
%8.280245416403167
%8.3453938090342478
%8.3072853725362129
%8.3772326468210268
%8.3021761619560794
%8.3973448224434441
%8.4257436538030372
%8.3105602650177267
%8.3266649773936844
%8.324988409318823
%8.269894527210702
% Connect4results.txt
Also, ${\sf agl}(\text{Connect 4})\approx 41.03(\pm 0.01)$ plies. 
%41.03479596645527
%41.02051504744945
%41.023884097185096
%41.024945262715093
%41.028995172932454
%41.030231147675153
%41.027407958464863
%41.034289903777122
%41.025213464069338
%41.027437516654238
%41.027165289211901
%41.015041702926048
While $P_1$ wins with \emph{perfect play}
(see \cite{Allis} and \url{tromp.github.io/c4/c4.html}), there is a caution:
it is likely that
 $P_1$ \emph{wins less overall}, at $\approx 27.71(\pm 0.21)\%$ 
%0.277268298653168
%0.27693688963681995
%0.27709621652567523
%0.27993785040030272
%0.27695953570791609
%0.27607684271822974
%0.27887482663733143
%0.2728586706318849
%0.27744802289299453
%0.27671794455372078
%0.27472547714904699
%0.28064601832621422
 than $P_2$ at $\approx 32.13(\pm 0.20)\%$ 
%0.32250165882075194
%0.3204563838677721
%0.32496473122275737
%0.3182108032259145
%0.31935644646851846
%0.32049581674427102
%0.32204756758579928
%0.31905422316109311
%0.31984228718016733
%0.32323199213943188
%0.3229372130714343
%0.32223387117376107
 (with draws at $\approx 40.16(\pm 0.30)\%$).
%0.40023004256420597
%0.40260672649540796
%0.39793905229030369
%0.40185134641378067
%0.40368401786293273
%0.40342734057771312
%0.39907760581521901
%0.40808710624604411
%0.40270968996703088
%0.40005006334602161
%0.40233730981840526
%0.39712011054047824

Finally, let ${\mathcal G}=$ \emph{Othello} (introduced into the United States in 1975 by Gabriel Industries, it is the modern version of \emph{Reversi}). This game is played on an $8\times 8$ board with disks $\dd$ and $\whd$, played by $P_1$ and $P_2$, respectively.
The rule is that $P_1$ places 
$\dd$ in a square if and only if there is another $\dd$ in the same row, column or diagonal
and $\whd$'s are contiguously between them. If the placement is valid, each of these  $\whd$'s flip to $\dd$'s. The same rule
applies to placing $\whd$ (except with the colors switched). A player may pass only if they do not have a move. The game ends when neither player has
a legal move. The winner is the one with the most disks. Finally, in Othello, the central squares start as $\tableau{\whd & \dd \\ \dd & \whd}$.

Na\"ively, ${\sf gts}({\text{Othello}})\leq 60!\approx 8.32\times 10^{81}$ (by filling all $60$ initially open squares in all possible ways). The {\sf gtc} estimate of
\cite{Allis} is $10^{58}$, based on an 
``in practice'' average game length of $58$ ply and an average of $10$ options/ply. 
Elsewhere, $10^{54}$ is estimated
for ${\sf gts}$, but without explanation/citation (see \url{en.wikipedia.org/wiki/Computer_Othello}).\footnote{\emph{A priori}, this is in contradiction with Allis' estimate, since by definition,
${\sf gtc}({\mathcal G})\leq {\sf gts}({\mathcal G})$. However, this can be reconciled as it looks like \cite{Allis} does not start with the 
four center squares filled.} 

One randomly generated Othello game $g$ ended with $\whd$ winning:
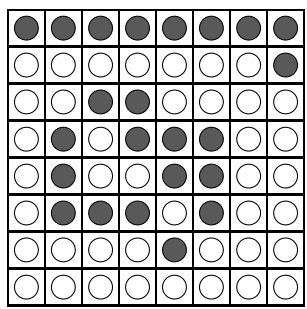
\[\tableau{\dd & \dd & \dd & \dd & \dd & \dd & \dd & \dd\\
\whd & \whd & \whd & \whd & \whd & \whd & \whd & \dd\\
\whd & \whd & \dd & \dd & \whd & \whd & \whd & \whd\\
\whd & \dd & \whd & \dd & \dd & \dd & \whd & \whd\\
\whd & \dd & \whd & \whd & \dd & \dd & \whd & \whd\\
\whd & \dd & \dd & \dd & \whd & \dd & \whd & \whd\\
\whd & \whd & \whd & \whd & \dd & \whd & \whd & \whd\\
\whd & \whd & \whd & \whd & \whd & \whd & \whd & \whd}\]
This gives $X(g)=2.49\times 10^{54}$. With $n=2\times 10^6$ ($24$ trials), (\ref{eqn:thetrick}) gives 
\[{\sf gts}({\text{Othello}}) \approx 6.47(\pm{0.19})\times 10^{54}.\]
%6.9012901434131941
%6.5464332421512057
%6.4354215061236435
%6.5173258839345713
%6.543463594060079
%6.2149756402498588
%6.237371074680193
%6.2930043663266665
%6.4474658169548495
%6.6521041918210933
%6.4947470753027611
%6.2637662478283232
%6.603929626806034
%6.2630889786806997
%6.300233628130404
%6.660524813472907
%6.365021523709854
%6.395288286557016
%6.389512047966175
%6.3780181391454148
%6.8813659438546075
%6.3147050248191067
%6.6260794274707196
%6.573058079465878
Also ${\sf agl}(\text{Othello})\approx 60.00(\pm 0.0004)$ ply, 
%59.997731040604478
%59.997792177117248
%59.998480183964709
%59.998522468422578
%59.998412762018276
%59.997778522527419
%59.998629586186262
%59.998761209880726
%59.998406567650676
%59.998570173710931
%59.99876049717124
%59.997882590199225
%59.998322635496187
%59.998171647960469
%59.998524213276255
%59.99787252632427
%59.99859030456465
%59.99833355598707
%59.99832121496595
%59.998576941394646
%59.998902964749291
%59.997067624762948
%59.998034457955782
%59.99788295265389
the draw rate is $\approx 4.95(\pm 0.30)\%$, 
%0.049006980126764514
%0.056356511400640463
%0.050499684529833555
%0.055617880765770365
%0.048428716977439563
%0.052653010268588235
%0.045451120953280966
%0.051825766027123674
%0.046618668077807422
%0.048288710869345566
%0.045245794747690561
%0.046798941256283449
%0.048598059973630169
%0.049412705423508514
%0.050738536384189438
%0.04849155593179061
%0.05227843881785069
%0.043293703201445204
%0.04839979780648971
%0.049981311391642927
%0.049839959407093612
%0.051384288874381262
%0.050503412736375437
%0.047766760794754184
$P_1$'s win rate is $\approx 43.36(\pm 1.56)\%$, 
%0.4445841408693636
%0.41451038241547455
%0.43197558199226249
%0.43146827567182361
%0.41807355745195857
%0.43250927107788067
%0.43706319757216783
%0.43709662828329859
%0.41883634336353703
%0.4517521236122633
%0.43745466324975263
%0.43566957353488212
%0.46833357151332894
%0.44148578116517317
%0.43227021131805432
%0.4322835400294229
%0.42317855550461814
%0.436376683945829
%0.4491738093270523
%0.44395518432347914
%0.3841296158029458
%0.43378887621596296
%0.43038912275725977
%0.44092457129329166
and $P_2$ has $\approx 51.69(\pm 1.51)\%$ of wins.
%0.50640887900389531
%0.52913310618398901
%0.51752473347791539
%0.51291384356238401
%0.53349772557066233
%0.51483771865353956
%0.51748568147454044
%0.5110776056895624
%0.53454498855863775
%0.49995916551838127
%0.51729954200251571
%0.51753148520885872
%0.48306836851305451
%0.50910151341134713
%0.51699125229773879
%0.5192249040387866
%0.5245430056775312
%0.5203296128527257
%0.502426392866458
%0.5060635042848779
%0.56603042478994581
%0.51482683490965553
%0.5191074645063648
%0.51130866791192
(Othello is currently unsolved, but the $4\times 4$ and $6\times 6$ versions have a forced win for $P_2$;
see \url{www.tothello.com}.)

Estimates for statistics of 
other games can be similarly attempted. Candidates include \emph{Checkers}, \emph{Dots and Boxes}, \emph{Go}, and \emph{Hex}. 
The second named author  has studied a simplified version of Chess, towards the 
understanding the difficulty of applying (\ref{eqn:thetrick}) to the full game. In addition, 
choosing each move \emph{uniformly} at random is inessential; one might wish to modify the probabilities, with the aim of reducing 
variance of the estimates. Finally, while antithetical to the crude approach espoused, if one does add evaluation and 
database information, (\ref{eqn:thetrick}) may be reinterpreted to estimate game tree \emph{complexity} and the number of ``sensible games''. Treatment of these topics will appear elsewhere.

\section*{Acknowledgements}
 
 We thank Yuguo Chen, Gidon Orelowitz and Anh Yong for helpful discussion.
AY was partially funded by a NSF grant and a Simons Collaboration Grant. This is a report for ICLUE, the Illinois Combinatorics Lab for Undergraduate Experience.

\end{document}